\documentclass{article}
\usepackage{latexsym}
\usepackage[]{amsfonts}
\usepackage{amsmath}
\usepackage{amssymb}

\title{\bf A Characterisation of $G_2(K)$}
\author{{\bf Christine Altseimer} \\ \\Institut f\"ur Mathematische
  Logik \\Albert-Ludwigs-Universit\"at
  Freiburg\\Eckerstr. 1\\ 79104 Freiburg\\ Germany}\date{} 

\begin{document} 
\pagestyle{plain}

\def\cc{\mathbb{C}}
\def\zz{\mathbb{Z}}
\def\ff{\mathbb{F}}
\def\nn{\mathbb{N}}
\def\qq{\mathbb{Q}}
\def\rr{\mathbb{R}}

\def\F{\mathop{\rm F}\nolimits}
\def\G{\mathop{\rm G}\nolimits}
\def\E{\mathop{\rm E}\nolimits}
\def\PSL{\mathop{\rm PSL}\nolimits}
\def\PSp{\mathop{\rm PSp}\nolimits}
\def\SL{\mathop{\rm SL}\nolimits}
\def\GL{\mathop{\rm GL}\nolimits}
\def\PGL{\mathop{\rm PGL}\nolimits}
\def\Sp{\mathop{\rm Sp}\nolimits}
\def\PSO{\mathop{\rm PSO}\nolimits}

\def\proof{ {\it Proof.} $\,$}

\newcommand{\semi}{\rtimes}  
\newcommand{\zs}{^{\circ}} 
\newcommand{\acf}{algebraically closed field}  
\newcommand{\fmr}{finite Morley rank} 
\newcommand{\uf}{if and only if}
\newcommand{\lemf}{\left(\begin{array}{cccc}} 
\newcommand{\lemt}{\left(\begin{array}{cc}}

\newcommand{\lemd}{\left(\begin{array}{ccc}}
\newcommand{\lemz}{\left(\begin{array}{cc}}
\newcommand{\rim}{\end{array}\right)}
\newcommand{\baz}{\begin{array}{cc}}
\newcommand{\bad}{\begin{array}{ccc}}
\newcommand{\ea}{\end{array}}

\newcommand{\ol}{\overline }
\newcommand{\lesm}{\left( \begin{smallmatrix}}
\newcommand{\rism}{ \end{smallmatrix} \right)}

\renewcommand{\baselinestretch}{1.1}

\newtheorem{lemma}{Lemma}
\newtheorem{theorem}[lemma]{Theorem}
\newtheorem{corollary}[lemma]{Corollary} 
\newtheorem{proposition}[lemma]{Proposition}
\newtheorem{remark}[lemma]{Remark}
\newtheorem{definition}[lemma]{Definition}
\newtheorem{fact}[lemma]{Fact}
\newtheorem{con}{Conjecture}

\renewcommand{\refname}{\centerline{\bf References}}

\maketitle

\begin{abstract} \thispagestyle{empty}
 This paper gives a characterization of the group $G_2(K)$ over
 some \acf\  $K$ of characteristic not 2 inside the class of simple
 $K^*$-groups of \fmr\ not interpreting a bad field using the
 structure of centralizers of involutions. This implies a general characterization of tame $K^*$-groups of odd type whose centralizers have a certain natural structure. 
\end{abstract}

\section{Introduction}
This paper belongs to a series of publications on the classification
of of tame simple groups of \fmr\ and odd type. The result and the
methods used to achieve it are intimatly related to
the characterisation of $PSp(4,K)$ given \cite{chr1}. The motivation
for this problem and all necessary definitions can be found there.   
The underlying conjecture is the following.

\begin{con}
\label{klitze}
Let $G$ be a simple tame $K^*$-group of odd type and Pr\"ufer 2-rank $2$.  
Let $i \in G$ be any involution and $C:=C_G(i)\zs/O(C_G(i))$. 
\begin{itemize}
\item[(1)] If\/ $C \cong \GL_2(K)$, then $G \cong \PSL_3(K)$
\item[(2)] If\/ $C \cong \PSL_2(K) \times K^* $, then $G \cong
  \PSp_4(K)$.  
\item[(3)] If\/ $C \cong \SL_2(K)*\SL_2(K)$, where the two copies of $\SL_2(K)$
  intersect non-trivially, then $G \cong \PSp_4(K)$ or $G\cong \G_2(K)$.   
\end{itemize}
Furthermore one of these three cases holds.
\end{con} 
 
We are going to prove cases $(2)$ and $(3)$ of the conjecture using \cite{chr1}. In \cite{diss} one can furthermore find a partial result of case $(1)$, if $O(C_G(i))=1$. To
prove the complete case $(1)$ seems to be extremely challenging and
one will probably need new methods. Especially we prove the following theorem:

\begin{theorem} 
\label{super}  
Let $G$ be a simple $K^*$-group of \fmr\ that does not
interpret a bad field. Assume that $G$ contains an involution $i$,
such that either     
\begin{itemize}   
\item[(i)] $C_G(i)\zs/O(C_G(i)) \cong \SL_2(K) * \SL_2(K)$ where the
  two copies of $\SL(2,K)$ intersect non-trivially or
\item[(ii)]$C_G(i)\zs/O(C_G(i)) \cong \PSL(2,K)\times K^*$ 
\end{itemize} 
for an \acf\ $K$ of characteristic neither $2$ nor $3$. Then $G \cong
\PSp_4(K)$ or $G \cong \G_2(K)$.   
\end{theorem} 

\section{Basic Results}

We are first going to show, that case $(ii)$ of Theorem \ref{super}
implies case $(i)$. 

\begin{proposition}
\label{restriction}
Let $G$ be a simple $K^*$-group of odd type that does not interpret any bad
field. Assume that $pr(G)=2$ and let $D$ be the four-subgroup which is
contained in the connected component of a Sylow 2-subgroup $S$ of
$G$. Set $O_D:=\bigcap_{l \in D^*} O(C_G(l))$. 
If\/ $T_0:=C_G(D)\zs/O_D \cong K^* \times K^*$ for
some \acf\ of characteristic not $2$, then $C_G(i)\zs/O(C_G(i))$ is
isomorphic to one of the following groups for any $i \in D^*$.
\begin{itemize}
\item[(i)] $K^* \times K^*$.
\item[(ii)] $\GL_2(K)$
\item[(iii)] $\PSL_2(K) \times K^* \cong \GL_2(K)/\langle -I \rangle$
  where $-I =   \left( \begin{smallmatrix} -1 & 0 \\ 0 & -1
    \end{smallmatrix}\right)$.  
\item[(iv)] $\SL_2(K)*\SL_2(K)$ where the two copies of $\SL_2(K)$
  intersect non-trivially.   
\end{itemize}
If\/ $C_G(i)\zs/O(C_G(i))$ is isomorphic to one of the groups in
$(iii)$ or $(iv)$, then $O(C_G(l))=1$ for all $l \in D^*$. 
\end{proposition}

\proof Let $i \in D^*$, then $C:=C_G(i)\zs/O(C_G(i))$ is a central
product of an abelian
divisible group $T$ and a semisimple group $H$ all of whose components
are simple algebraic groups over algebraically 
closed fields of characteristic different from 2 by \cite{nato}. Furthermore $pr(C)=2$ by
\cite{chr2} and $pr(H)=pr(C)-pr(T)$ by Lemma
\cite{chr2}. Thus $H$ is the central product of simple algebraic groups
of Pr\"ufer 2-rank less than or equal to 2 by
\cite[27.5]{hum75}. 

Set $C_k:=C_G(k)$ and $O_k:=O(C_k)$ for all $k
\in I(G)$. Let $\theta$ be the signalizer functor defined by
$\theta(s):=O_s$ for all $s \in I(G)$. Then $O_i \cap C_j = O_j
\cap C_i = O_j \cap O_i = O_D$ for all $j \in D^*$, $j \neq i$. Hence 
\[C_i \cap C_j / O_D \cong (C_i \cap C_j)O_i/O_i \cong (C_i \cap
C_j)O_j/O_j.\] 
Thus by  \cite[2.52]{mixed} $(C_i \cap C_j) / O_D \cong C_{C_i/O_i}(jO_i)
\cong C_{C_j/O_j}(iO_j)$. This implies that
 $T \leq C_{C}(jO_i) \cong C_G(D) / O_D = T_0$ and either $T =1$, $T
 \cong K^*$ or $T  \cong T_0$, since $T$ is connected and $G$ does not
 interpret a bad field.  

Thus $(i)$ holds, if $pr(T)=2$, as in this case $H=1$ and $T \cong
T_0$. If $pr(T)=1$, then $T \cong K^*$ and $pr(H)=1$. Hence $H$ is of
type $A_1(K)$, i.e.\/ $H \cong \SL_2(K)$ or $H \cong \PSL_2(K)$. 
If $C=H \times T$, then $H \cong \PSL_2(K)$, as $C_{C}(jO_i) \cong
T_0$ for any involution $j \in H$. Thus $C \cong\GL_2(K)/\langle -I \rangle$
in this case. If on the other hand $C = H * T$, where $H \cap T$ is
non-trivial, then $H \cong \SL_2(K)$ and $C \cong \GL_2(K)$.

If finally $pr(H)=2$, then $T$ is
trivial by \cite{chr2} and $H$ is of type  $A_1(K)
\times A_1(K)$, $A_2(K)$, $C_2(K)$ or $\G_2(K)$. However, $H$ has to
contain a central involution. Thus $H$ cannot be of type $A_2(K)$ or
$\G_2(K)$ and it cannot be isomorphic to $\PSp_4(K)$. As furthermore
$C_{C}(jO_i) \cong T_0$ for any $j \in D\backslash \langle i
\rangle$, $H$ cannot be isomorphic to $\Sp_4(K)$ either. Hence $H$ is
of type $A_1(K) \times A_1(K)$, contains the central involution $i$
and $C_{C}(jO_i) \cong T_0$ for any other involution $j$. Especially
$C \cong \SL_2(K) * \SL_2(K)$, where the two copies of $\SL_2(K)$
intersect in $\langle i \rangle$. 

If $C$ is isomorphic to the group in $(iii)$ or $(iv)$, then
$C_G(i)\zs$ contains an elementary abelian subgroup $E$ of order $8$
which contains $D$ by \cite{chr1} as $C$ does. Thus
$O(C_G(i))=1$ by \cite[17]{chr2}. \hfill $\Box$  

\begin{lemma}
\label{ratherin}
Let $G$ be a $K^*$-group of \fmr. Assume that there exist two definable
subgroups $N, H \leq G$ such that $N \lhd H$ and $\ol H:=H/N$  is a proper
simple section of $G$ allowing no graph automorphisms. Then $N_G(H) =
C_{N_G(H)}(\ol H)H$. 
\end{lemma} 

\proof  Let $g \in N_G(H)$ and set $R:=d(g)$. $R \leq N_G(H)$ acts on
$\ol H$. As $G$ is a $K^*$-group and $H$ is a proper simple section,
$H$ is a simple algebraic group over an \acf. Consider
the semidirect product $\ol H \semi R/C_R(\ol H)$. Then, viewing $R/C_R(\ol
H)$ as a subgroup of $Aut(\ol H)$,  $R/C_R(\ol H) \leq Inn(\ol H) \Gamma$,
where $\Gamma$ are the graph automorphisms of $\ol H$ by
\cite[8.4]{BN}. As $\Gamma =1$,  either $g \in C_{N_G(H)}(\ol H)$ or
there exists an $h \in  H$ such that $x^g \in x^hN$
for all $x \in H$ and $gh^{-1} \in C_{N_G(H)}(\ol H)$. Thus $g \in
C_{N_G(H)}(\ol H)H$. 
\hfill $\Box$ \vspace{1ex}
 
\begin{lemma}
\label{cicacen}
Let $H$ be a group of \fmr, such that $H\zs$ is a reductive algebraic group and
let $\ol H\zs:=H\zs/Z(H\zs)$. Then $C_{C_H(Z(H\zs))}(H\zs) =
C_{C_H(Z(H\zs))}(\ol H\zs)$. 
\end{lemma}

\proof  Let $g \in C_H(\ol H\zs) \cap C_H(Z(H\zs))$.  As $H\zs$ is reductive,
$\ol H\zs$ is a simple group by \cite[27.5]{hum75} and $m:=o(\ol g) <
|H/H\zs|$. Set
$z_x:=x^{-1}x^g \in Z(H\zs)$ for all $x \in H\zs$. Then, as
$x^g=xz_x$, $x=x^{g^m}= xz_x^m$ and $o(z_x)|m$ for all $x \in
H\zs$. On the other hand $(x^m)^g=(xz_x)^m = x^m z_x^m$ and $g \in
C_H(x^m)$ for all $x \in H\zs$.   Especially $C_P(g)=P$, for any
maximal torus $P$ of $H\zs$, since maximal tori are divisible. As finally
$H\zs$ is a reductive algebraic group, $H \zs = \langle P^h | \; h \in
H\zs \rangle \leq C_H(g)$ by \cite[ex.\ 12, p.\ 162]{hum75} and $g \in
C_H(H\zs)$. \hfill $\Box$ 
 
\begin{theorem}
\label{exclude}
Let $G$ be a simple $K^*$-group of \fmr\ that does not
interpret a bad field. Assume that $C_G(i)\zs/O(C_G(i)) \cong \PSL_2(K) \times K^*$ for an involution $i \in G$ where $K$ is an \acf\ 
of characteristic not $2$. Then $G$ contains an involution $j$, such
that $C_G(j)\zs \cong \SL_2(K) * \SL_2(K)$.
\end{theorem}

\proof Let $S$ be Sylow 2-subgroup of
$G$ that contains $i$ and $D \leq S\zs$ a four-subgroup. We may assume
that $i \in D$ by Proposition \ref{restriction}. Set $T:=C_G(D)\zs$.
Hence $O(C_G(t))=1$ for all $t \in D^*$ by Proposition
\ref{restriction} again.  Furthermore $G=\langle C_G(k)\zs| \; k \in
D_1^*\rangle$ for any four-subgroup $D_1 \leq G$ by \cite{chr2}. Set $T:=C_G(D)\zs$ and let $u \in (C_G(D) \cap
C_G(i)\zs) \backslash T$. Then $C_T(u)=Z(C_G(i)\zs)$.

Assume that $C_G(j)\zs \cong \GL_2(K)$ for some $j \in D^*$.
Then there exists an involution
$w \in C_G(j)\zs \cap N_G(D)$ such that $i^w =ij$. Then
$uu^w=w^uw \in C_G(D) \cap C_G(j)\zs =T$. Furthermore, as $u$
centralises $Z(C_G(i)\zs)$, $u^w$ centralises $Z(C_G(i)\zs)$ as
well. On the other hand
$C_T(u^w) = C_T(u)^w = Z(C_G(ij)\zs)$. Contradiction as $Z(C_G(ij)\zs)
\neq Z(C_G(i)\zs)$. 

Assume that $C_G(j)\zs \cong \PSL_2(K)\times K^*$ for some $j
\in D^*$, $j \neq i$. Let $v \in
C_{C_G(j)\zs}(D)\backslash T$ be an involution. As
$C_T(v)\zs=Z(C_G(j)\zs)$, $v \in C_G(i) \backslash C_G(i)\zs$.  Let
$C:=C_G(i)\zs$ and $\ol C:=C/Z(C)$. By Lemma \ref{ratherin}
$C_G(i)=C_{C_G(i)}(\ol C)C$. 

$v$ normalises $Z(C) \cong K^*$ and $G$ does not interpret a bad
field. Thus $v$ either inverts $Z(C)$ or centralizes it by
\cite[10.5]{BN}. Since $C_T(v)=Z(C_G(j)\zs)$, the second case cannot occur and
$v$ inverts
$Z(C)$. Let $x \in C$, such that $xv \in C_{C_G(i)}(\ol C)$. Then $x
\in N_G(D)\cap C= T \langle u \rangle$ and $xv$ inverts $Z(C)$. As $u$
inverts $Z(C_G(j)\zs)$ as above, $yv$ inverts $T$ for all $y \in
uT$. Hence $x \in T$, since $xv \in C_{C_G(i)}(\ol C)$. Furthermore we
may assume that $x
\in Z(C_G(j)\zs)$, as $T=Z(C_G(j)\zs)Z(C)$. Thus $(xv)^2=x^2 \in
C_{C}(\ol C)=Z(C) \cap Z(C_G(j)\zs) =1$ by Lemma \ref{cicacen} and $x
\in \langle j \rangle$. Then $[xv,C]
\subseteq Z(C)$ and $C=Z(C)*C_C(xv)$ by \cite[ex.\ 10, p.\ 98]{BN}
where $Z(C) \cap C_C(xv)= \langle i \rangle$. Set 
$H:=C_C(xv)\zs$. As $C \cong \PSL_2(K)\times K^*$, $H \cong
\PSL_2(K)$ and $C=Z(C) \times H$. As $x
\in \langle j \rangle$, $xv \in I(C_G(D))$. Thus $D_1:=\langle i, xv
\rangle$ is a four-subgroup and $G:=\langle C_G(k)|\; k \in D_1^*
\rangle$. Furthermore $H \leq C_G(D_1)$. As $H$ is a normal subgroup
of $C_G(i)\zs$ and $G$ is simple, $H$ cannot be a normal subgroup of
both $C_G(ux)$ and $C_G(uxi)$. As $ux$ and $uxi$ are conjugate to
involutions in $D$ by elementary computation in $\GL_2(K)/<-I>$, this implies that $C_G(k) \cong \SL(2,K)*\SL(2,K)$ for $k=ux$ or $k=uxi$ by Proposition
\ref{restriction}.  

As finally $G=\langle C_G(k)\zs | \; k \in D^* \rangle$, it is impossible that
$C_G(k)\zs$ is abelian for all $k \in D \backslash \langle i \rangle$. 
This implies the claim by Proposition \ref{restriction}. \hfill $\Box$ \vspace{1ex}

Thus Theorem \ref{super} basically consists of two parts. The first one is the characterization of $\PSp_4$ as in \cite{chr1}. 

\begin{fact} \index{PSP@$\PSp_4(K)$} 
\label{main} 
Let $C_1$ and $C_2$ be the non-isomorphic
  centralizers in $\PSp_4(K)$ of two involutions where K is an \acf\
  of characteristic not $2$.  Let $G$ be
a $K^*$-group of \fmr\ that does not interpret a bad field. Assume
that $G$ contains two involutions $i$ and $j$ such that
$C_G(i)\cong C_1$ and $C_G(j) \cong 
C_2$.  Then $G \cong \PSp_4(K)$. 
 \end{fact} 

The corresponding result for $\G_2(K)$ is the following theorem, which
will be proven in the following section.

\begin{theorem}   \index{G2@$\G_2(K)$}
\label{maing2} 
Let $C$ be the centralizer in $\G_2(K)$
  of an involutions where K is an \acf\   of characteristic not $2$. Let $G$ be a
$K^*$-group of \fmr\ that does not interpret a bad
  field. Assume that $G$ contains one conjugacy class $i^G$ of involutions
such that $C_G(i)\cong C$.  Then $G \cong \G_2(K)$ if $char(K) \neq 3$.
 \end{theorem}

\section{A characterization of $\G_2(K)$}

In this section we prove Theorem \ref{maing2}. Let $G$ be as in
Theorem \ref{maing2}, $D:=\langle i_0, i_1 \rangle$ a four-subgroup which is
contained in the connected component of a Sylow 2-subgroup of $G$ and set
$i_2=i_0i_1$. Let furthermore $T:=C_G(D)\zs$. $C_G(i_0)$ contains
four quasiunipotent subgroups which are normalized by $T$ and 
isomorphic to $K^+$. Let $X$ and $Y$ be two of them that centralize
each other and let $v$ be an involution in $N_{C_G(i_0)}(T) \backslash
T$ that normalizes $X$. Write $\ol x:=xT$ for all $x \in N_G(T)$.

\begin{lemma}
\label{weyl}
Let $W:=N_G(T)/T$. Then $|W|=12$ and $W$ is generated by two
involutions $\ol w$ and $\ol v$, where $w \in I(C_G(i_1))$, such that
$i_0^{wv}=i_1=i_2^{vw}$ and $o(wv)=6$. The 
center of $W$ is $\langle \ol z \rangle$ where $z:=(wv)^3$. Set
$y:=(wv)^4=(wv)z$ and $w_{k+1}:=w^{y^k}$, $v_k:=v^{y^k}$ for $k \in
\nn$. Then $\ol w_k = \ol {w_{k+3}}$ and $\ol v_{k} = \ol{v_{k+3}}$. The
elements in $W$ are hence 
\[\ol 1, \ol z = \ol{(wv)^3}, \ol{wv}=\ol{yz},
\ol{(wv)^2} = \ol{y^{-1}}, \ol{(wv)^4}=\ol y,
\ol{(wv)^5} = \ol{vw} = \ol{y^{-1}z}\] 
and 
\[\ol w_{k+1}=\ol{w^{(vw)^{2k}}}=\ol{y^{k}w} \mbox{ and } \ol v_k =
\ol{v^{(vw)^{2k}}} = \ol{y^{k}v} \] 
for $k=0,1,2$.
\end{lemma}   

\proof As $G$ contains one conjugacy class of involutions, all
involutions of $D$ are conjugate in $N_G(T)$ by
\cite[10.22]{BN}. $N_G(D)/C_G(D) \cong S_3$, $N_G(T)=N_G(D)$ and
$|C_G(D)/C_G(T)| = 2$, which implies that $W$ is the dihedral group of
order 12.  \hfill $\Box$

\begin{lemma} 
\label{grundbn3}
Let $w_0, v_0$ as in Lemma \ref{weyl}. Then $w_0TXw_0 \subseteq
TXw_0X$ and $v_0TYv_0 \subseteq TYv_0Y$. Furthermore $C_G(i_0)\zs =
\langle T, X, Y ,v_0, w_0 \rangle$.  
\end{lemma}

\proof  $L:=\langle X, T, X^{w_0}\rangle$ is a reductive algebraic
group of rank $1$. $L$ has thus a $BN$-pair $(B_1,N_1)$, where $B_1 := TX$ and
$N_1= \langle w_0, T \rangle$. Hence $w_0TXw_0 \subseteq
 TX \cup TXw_0X = TXw_0X$ and $L:=\langle T, X, w_o \rangle$. The
 same argument for $\langle Y ,T ,Y_0 \rangle$ yields the remaining part of
 the lemma. \hfill $\Box$ 

\begin{lemma}
\label{xconju}
Let $X_1:=X$, $Y_1:=Y$, $X_2:=X_1^{w_0}$, $Y_2:=Y_1^{v_0}$ and
$X_n^{(\lambda)}:= X_n^{y^{\lambda}}$ as well as $Y_n^{(\lambda)}:=
Y_n^{y^{\lambda}}$ for any $n=1,2$ and $\lambda \in \zz$. Then 
\begin{itemize}
\item[(i)] $X_n^{(\lambda)}=X_n^{(\lambda+3)}$ and
 $Y_n^{(\lambda)}=Y_n^{(\lambda+3)}$ for all $n=1,2$ and
 $\lambda \in \zz$. 
 \item[(ii)] $v_{\lambda}$ centralizes $X_n^{(\lambda)}$ and
 $w_{\lambda}$ centralizes $Y_n^{(\lambda)}$  for $\lambda = 0,1,2$
 and $n=1,2$.
\item[(iii)] $ (X_1^{(\kappa)})^{w_{\lambda}} =  X_2^{(2\kappa -
 \lambda)}$ and $(X_1^{(\kappa)})^{v_{\lambda}} =  X_1^{(2\kappa -
 \lambda)}$ for $\kappa, \lambda = 0, 1, 2$.
 \item[(iv)] $ (Y_1^{(\kappa)})^{w_{\lambda}} =  Y_1^{(2\kappa -
 \lambda)}$ and $ (Y_1^{(\kappa)})^{v_{\lambda}} =  Y_2^{(2\kappa -
 \lambda)}$ for $\kappa, \lambda = 0, 1, 2$.
\end{itemize}
\end{lemma}

\proof Since $\ol{y^3} = \ol 1$, $(i)$ follows. We have furthermore
 chosen $w_0$ such that $w_0$ centralizes $Y_1$. Hence $v_0$ has to
 centralize $X_1$ which yields $(ii)$. To prove $(iii)$ note that
 $\ol{ y^{\lambda}}=\ol{w_\lambda w_0}=\ol{v_{\lambda}v_0}$ and hence 
\[\ol{y^{\kappa}w_{\lambda}y^{\lambda - 2 \kappa}w_0} = \ol{y^{3
 \kappa - \lambda}w_{\lambda}w_0}= \ol{y^{-\lambda}y^{\lambda}}= \ol
 1\]
and
\[\ol{y^{\kappa}v_{\lambda}y^{\lambda - 2 \kappa}v_0} = 
\ol{y^{3\kappa -\lambda}v_{\lambda}v_0}= \ol 1\]
for $\lambda, \kappa \in 0,1,2$. Since however $X_1^{w_0}=X_2$,
$X_1^{v_0}=X_1$, $Y_1^{w_0}=Y_1$ and $Y_1^{v_0}=Y_2$ by $(ii)$,
$(iii)$ and
$(iv)$ follow. \hfill $\Box$

\begin{lemma}
\label{finitealg}
Let $G$ be a connected $K$-group of \fmr\ such that the solvable
radical $\sigma$ of $G$ is finite. Then $G$ is a central product of
quasi-simple algebraic groups over \acf s. 
\end{lemma}

\proof Since any definable action of a definable connected group on a
finite set is trivial $\sigma = Z(G)$. Furthermore $G/\sigma$ is
isomorphic to the direct product of simple algebraic groups over
algebraically closed fields by \cite{alt94}. Then
$G/\sigma=(G/\sigma)'=G'\sigma/\sigma \cong G'/(G' \cap \sigma)$ and
$rk(G)=rk(G')$. As $G$ is connected, $G=G'$ and $G$ is semisimple.
Assume that $G$ is quasi-simple. Then $G$ is an algebraic group by
\cite{central}. Let now $G/\sigma \cong A_1 \times \cdots \times A_m$
for some $m \in \nn$, where $A_i$ are simple algebraic groups
over \acf s for $1 \leq i \leq m$. Let $G_i$ be the preimages of $A_i$
in $G$ for $1 \leq i \leq m$. Then $G=G_1\zs....G_m\zs$. Now $[G_i\zs
,G_j \zs]$ is a connected subgroup of $\sigma$ and hence trivial for
any $1 \leq i < j \leq m$. Furthermore $[G_i\zs,G_i\zs]=G_i\zs$ for
any $1 \leq i \leq m$ as above. Thus $G$ is a central product of $m$
quasi-simple algebraic groups by the first case. \hfill $\Box$ 

\begin{proposition}
\label{subgroups}
  If\/ $char(K) \neq 3$ then there exists an element $t$ of order 3, such
  that  $C_G(t)\zs \cong \SL_3(K)$. Actually $C_G(t)\zs \cong \SL_3(K)$
  for exactly two elements of order 3 in $T$. 
\end{proposition}

\proof Assume that $char(K)$ is not 3. Since $y^3 \in T$, $y=y't$, where we can
choose $y'$ to be a 3-element by \cite[ex.\ 11, p.\ 93]{BN}. Let $P$
be a Sylow 3-subgroup of $N_G(T)$.  Since $N_G(T)$ is solvable, $P$ is
nilpotent-by-finite by \cite[6.20]{BN}. Hence there exists an element
$t \in T \cap Z(P)$ by \cite[ex.\ 12, p.\ 14]{BN}. It follows that $t$
is centralized by  $y$. Thus $(t^{w_0}t)^{y^2} = (t^{w_0}t)^{w_2w_0} =
t^{yw_0}t =t^{w_0}t$ and $t^{w_0}t \in C_T(y)$. Hence either
$t^{w_0}=t^{-1}$ or $t^{w_0}t \in C_T(y)$ is an element of order 3 that is
inverted by $w_0$. We may assume that $t$ is inverted by $w_0$. Now
$C_G(t)\zs = \langle C_{C_G(t)}(l)\zs |\; l \in D^* \rangle$ and
$C_{C_G(i_0)}(t)=C_{C_G(i_1)}(t)^{y^2 }=C_{C_G(i_2)}(t)^{y}$.
However, $C_{C_G(i_0)}(t)\zs  \cong \GL_2(K)$. 

We show that the solvable radical $\sigma$ of $C_G(t)\zs$ is
finite. $\sigma$ is normalized by $T$ and hence $\sigma\zs = \langle
C_{\sigma}(l)\zs | l \in D \rangle$ by \cite[4.12]{nato}. However,
$C_{\sigma}(l)\zs \leq Z(C_{C_G(l)}(t)\zs)$ for all $l \in D^*$ as
$C_{\sigma}(l)\zs$ is contained in the solvable radical of
$C_{C_G(l)}(t)\zs$. Thus $\sigma\zs \leq \bigcap_{l \in
  D^*}Z(C_{C_G(l)}(t)\zs) \leq T$. One the other hand
$Z(C_{C_G(i_0)}(t)\zs) \cong K^*$ and as $G$ does not interpret a bad
field, $\sigma$ is either finite or $i_0 \in
Z(C_{C_G(i_0)}(t)\zs)=Z(C_{C_G(i_1)}(t)\zs)$. The second case cannot
occur since $C_G(D)\zs=T$ and $\sigma$ finite. 

Thus $C_G(t)\zs$ is a central product of quasi-simple algebraic groups
by Lemma \ref{finitealg} and $C_G(t)\zs \cong \SL_3(K)$. Furthermore
$T \leq C_G(t) \zs$ contains eight elements of order 3 and while $y$
centralizes exactly $t$ and
$t^{-1}$, $yw$ operates transitively on the remaining six. \hfill $\Box$

\noindent\\ The following two propositions can be proven exactly as the
corresponding results in \cite{chr1}.

\begin{proposition} 
\label{maxuni} 
There exists a subgroup $V \leq G$ such that $V$ is a maximal
quasiunipotent group and is normalized by $T$.   
\end{proposition}

\begin{proposition} 
\label{NU} 
Let $V$ a maximal quasiunipotent group
of $G$ which is normalized by $T$.  Then $N_G(V)\zs = V \semi T$.  
\end{proposition}

\noindent\\ Now we can construct a $BN$-pair as in \cite{rusbor}

\begin{proposition}
\label{thisis}
If\/ $char(K) \neq 3$ then there exists a maximal quasiunipotent
group $Q$ in $G$ which is normalized by $T$ such that 
\begin{itemize}
\item[(i)] $Q=V \semi Y_1$ where $V$ is normalized by $v_0$ and
\item[(ii)] $Q=M \semi X_2^{(1)}$ where $M$ is normalized by $w_1$.
\end{itemize}
Furthermore $Q=\langle X_1, X_2^{(1)}, X_2^{(2)}, Y_1, Y_2^{(1)},
Y_1^{(2)} \rangle$.
\end{proposition}

\proof Since $char(K)$ is not $3$, there exists an element $t \in T$ of
order 3, such that $C_G(t)\zs \cong \SL_3(K)$. We can hence assume
that $C_G(t)\zs$ contains the maximal quasiunipotent subgroup
$U:=X_1\langle X_2^{(1)}, X_2^{(2)}\rangle$, where $X_1$ is central in
$U$. Consider $C_G(X_1)$. We proceed as in \cite{chr1}. $U \leq C_G(X_1)$ and $L:=\langle T, Y_1, Y_2 \rangle \leq
C_G(X_1)$. Thus $C_G(X_1)\zs$ is not solvable. Furthermore
$C_G(X_1)\zs$ has Pr\"ufer 2-rank 1, since $C_G(T)\zs =T$. Thus
$C_G(X_1)\zs/\sigma \cong \PSL_2(K)$ by \cite{alt94} and
\cite{chr2}, where $\sigma$ is the solvable radical of
$C_G(X_1)$. On the other hand $L\sigma/\sigma  \cong L/(\sigma \cap L)
\cong L/Z(L) \cong \PSL_2(K)$. Hence $C_G(X_1)\zs =L\sigma\zs$ and
$U\leq \sigma$. Set $V:=Q(\sigma)$. Then $V^{v_0}=V$. 

Assume that $U=V$. Then $Y_1U$ is a quasiunipotent group by \cite{qua}. Especially $Y_1U$ is nilpotent and $W:=Y_1X_1 < Y_1U$ has
infinite index in $N:=N_{Y_1U}(W)\zs$ by \cite[6.3]{BN}. $N$
is a quasiunipotent group that is normalized by $T$ and hence $N =
\langle C_N(l)\zs |\; l \in D^*\rangle$ by
\cite[4.6]{nato}. Furthermore $C_N(i_0)=W$ and $C_N(i_k)$ is
either trivial or equals $X_2^{(k)}$ for $k=1,2$ as $T$ normalises $N$
and acts transitively on $X_2\backslash \{1\}$. We may assume
that $X_2^{(1)} \leq N$. 

Set $P:=WX_1^{(2)} \semi T$. Then $P \leq N_G(W)$. Let $w \in W
\backslash\{1\}$ such that $C_G(w) \cap T = \langle i_0 \rangle$. As
$T$ acts transitively on $W \backslash \{1\}$,
for any element $x \in P$, there exists an element $t_x \in T$ such
that $w^x=w^{t_x}$. Thus $P \subseteq C_P(w)T$. As $i_1$ and $i_2$ invert
$W$, $C_P(w)$ is a solvable group normalized by $D$. Thus
$C_P(w)\zs=\langle C_{C_P(w)}(l)\zs |\; l \in D^* \rangle$ by
\cite[4.6]{BN}. Furthermore $(C_P(w) \cap C_G(i_0))\zs = W$, $(C_P(w)
\cap C_G(i_1))\zs \leq X_1^{2}$ and $(C_P(w) \cap C_G(i_2))\zs =
1$. Thus $Q:=C_P(w)\zs \leq WX_1^{(2)}$ is a quasiunipotent subgroup of
$P$ containing $W$. Since there is a definable surjective map from
$C_P(w) \times T$ onto $N_P(W)$ 
\[rk\big(WX_1^{(2)}\big)+ rk(T) = rk\big(N_P(W)\big) \leq
rk\big(C_P(w)\big)+rk(T) = rk(Q) + rk(T).\]
Thus $Q=WX_1^{(2)}$ and $N_P(W)=C_P(w)\zs \semi T$. Especially $Q \leq
C_P(w^t)$ for all $t \in T$ and hence $Q \leq C_G(W)$ as $T$ acts
transitively on $W$. Thus $X_2^{(1)} \leq C_G(W)$ and
$X_2^{(1)}$ centralizes $Y_1$. This implies that $Y_1^{(1)} = Y_1^{w_2} \leq
C_G(X_2^{(1)})^{w_2}=C_G(X_1)$. As $C_G(X_1)\zs =L\sigma$ this implies
that $Y_1^{(1)} \leq V$. Contradiction. 

Hence $U<V$ and either $\langle
Y_2^{(1)}, (Y_2^{(1)})^{v_0}\rangle \leq V$ or  $\langle
Y_1^{(1)}, (Y_1^{(1)})^{v_0}\rangle \leq V$. We may assume that $V= \langle
Y_2^{(1)}, Y_1^{(2)}, U\rangle$. Set $Q:=Y_1V$. $Q$ is
obviously a maximal quasiunipotent subgroup of $G$ and $V$ is
normalized by $v_0$ by construction. Set $Y:=\langle Y_1,
Y_2^{(1)}, Y_1^{(2)}\rangle$. Then $Y$ is invariant
under $w_1$ and $M:=Q \cap Q^{w_1} = \langle Y
, X_1 X_2^{(2)}\rangle$. Furthermore $M$ is
normalized by $X_2^{(1)}$, since $N_Q(M)\zs=Q$ by \cite[6.3]{BN} as in
the previous paragraph. \hfill $\Box$

\begin{theorem} 
\label{char3} 
If\/ $char(K) \neq 3$, $G$ is a split BN-pair of Tits rank 2, where
$B:=N_G(Q)\zs$, $N:=N_G(T)$ and $\overline w_2, \overline v_0$ are the
generators of the Weyl group.   
\end{theorem}

\noindent This proves our result, Theorem \ref{maing2}, by
\cite{qua}, noting that the gap in the proof of the necessary theorem is filled by \cite{kramer}.\vspace{1ex} 

\proof $G=\langle C_G(k)|\; k \in D^*
\rangle$ by \cite[18]{chr2} and $\langle
N,X_1Y_1\rangle \geq C_G(i_0)$ by Lemma \ref{grundbn3}. Thus
$G=\langle B, N \rangle$ as $i_1, i_2$ are conjugate to $i_0$ in $N$. BN1 now
follows since by Proposition \ref{NU} $B=QT$ and thus 
$B\cap N = T \lhd N$. Furthermore $\langle \overline v_0, \overline
w_1 \rangle  = N_G(T)/T$ by Lemma \ref{weyl} proving BN2. Let $S:=\{
\overline v_0, \overline w_1\}$ and $W:=N/T$.  

To prove BN3 we will show that $vBw \subseteq BvwB \cup BvB$ for
all $v,w \in N$ such that $\overline v \in W$ and $\overline w
\in S$. By Lemma \ref{xconju} and Proposition
\ref{thisis}  
\[1Bv_0 \subseteq B1v_0B  \]
\[v_0Bv_0 = v_0TY_1Vv_0  \subseteq TY_1v_0Y_1V \subseteq Bv_0B \] 
\[v_2Bv_0 = v_2TY_1Vv_0 = TY_2^{(1)}v_2v_0V \subseteq Bv_2v_0B  \]
\[w_0Bv_0 = w_0TY_1Vv_0 = TY_1w_0v_0V \subseteq Bw_0v_0B\]
\[w_1Bv_0 = w_1TY_1Vv_0 = TY_1^{(2)}w_1v_0V \subseteq Bw_1v_0B  \]
\[y^2Bv_0 = y^2TY_1Vv_0 = TY_1^{(2)}y^2v_0V \subseteq By^2v_0B  \]
\[v_0w_1Bv_0 = v_0w_1TY_1Vv_0 = TY_2^{(1)}v_0w_1v_0\subseteq
Bv_0w_1v_0B  \]
and thus
\[zBv_0 = (zv_0)v_0Bv_0 \subseteq w_0Bv_0B \subseteq Bw_0v_0B = BzB
\]
\[w_1v_0Bv_0 \subseteq w_1Bv_0B \subseteq Bw_1v_0B  \]
\[w_2Bv_0=(w_2v_0)v_0Bv_0 \subseteq v_0w_1Bv_0B \subseteq
Bv_0w_1v_0B
= Bw_2B.\]
Finally
\[yBv_0 = yTY_1Vv_0 = TY_1^{(1)}yv_0V = TY_1^{(1)}v_1V =
Tv_1Y_1^{(2)}V  \subseteq Byv_0B \]
which gives us
\[v_1Bv_0 = (v_1v_0)v_0Bv_0 \subseteq yBv_0B \subseteq
Byv_0B=Bv_1B.\]

\noindent On the other hand
\[1Bw_1 \subseteq B1w_1B  \]
\[w_1Bw_1 = w_1TX_2^{(1)}Mw_1  \subseteq TX_2^{(1)}w_1X_2^{(1)}M
\subseteq Bw_1B \]  
\[v_1Bw_1 = v_1TX_2^{(1)}Mw_1 = TX_2^{(1)}v_1w_1M \subseteq
Bv_1w_1B  \]
\[w_2Bw_1 = w_2TX_2^{(1)}Mw_1 = TX_1w_2w_1M \subseteq
Bw_2w_1B\]
\[v_0Bw_1 = v_0TX_2^{(1)}Mw_1 = TX_2^{(2)}v_0w_1M \subseteq
Bv_0w_1B  \]
\[yBw_1 = yTX_2^{(1)}Mw_1 = TX_2^{(2)}yw_1M \subseteq Byw_1B  \]
\[w_1v_0Bw_1 = w_1v_0TX_2^{(1)}Mw_1 = TX_1w_1v_0w_1\subseteq
Bw_1v_0w_1B  \]
and thus
\[zBw_1 = (zw_1)w_1Bw_1 \subseteq v_1Bw_1B \subseteq Bv_1w_1B =
BzB\]
\[v_0w_1Bw_1 \subseteq v_0Bw_1B \subseteq Bv_0w_1B  \]
\[w_2Bw_1=(w_2w_1)w_1Bw_1 \subseteq yBw_1B \subseteq
Byw_1B=Bw_2B. \]
Finally
\[y^2Bw_1 = y^2TX_2^{(1)}Mw_1 = TX_2y^2w_1M = TX_2w_0M =
Tw_0X_1M  \subseteq By^2w_1B \] 
which gives us
\[v_2Bw_1 = (v_2w_1)w_1Bw_1 \subseteq w_1v_0Bw_1B \subseteq
Bw_1v_0w_1B=Bv_2B.\]

\noindent and BN3 is valid. By Proposition \ref{thisis} again
$Q^{v_0}=VY_2 \neq VY_1 = Q$ and $Q^{w_2}=MX_1^{(1)} \neq
 MX_1^{(2)}=Q$ which proves BN4 and gives us the theorem. \hfill  $\Box$ 
 
 \index{group!of finite Morley rank|)}   \index{G2@$\G_2(K)$|)}

\section{Proof of Theorem \ref{super}}

This section is devoted to the proof that Theorems
\ref{main} and \ref{maing2} imply Theorem \ref{super}. \vspace{1ex}

Let $G$ be a group as in Theorem \ref{super}. By Theorem \ref{exclude}
there exists an involution $i \in G$, such that $C_G(i)\zs/O(C_G(i)) \cong L_1 * L_2$ where $L_n \cong \SL_2(K)$ for an \acf\ $K$ of characteristic
neither $2$ nor $3$ and $n=1,2$. We may assume that $i \in C_G(i)\zs$ by Lemma
\ref{restriction}. Thus $L_1 \cap L_2 = \langle \ol{\mbox{\it{\i}}}
\rangle$. Let $S$ be a Sylow 2-subgroup of $C_G(i)$, $D$ a
four-subgroup such that $D \leq S\zs$ and set $T:=C_G(D)\zs$. Then $D$
contains a central involution of $S$ by \cite[ex.\ 12, p.\ 14]{BN} which means that $i$ is central in $S$, since $j$ and $ij$ are conjugate in $S$
for $j\in D\backslash \langle i \rangle$.  

\begin{proposition}
\label{corset}
$G=\langle C_G(k)\zs| \; k \in D_1^*\rangle$ for any four-subgroup
$D_1 \leq G$ and $O(C_G(k))=1$ for all $k \in D^*$.
\end{proposition}

\proof $O(C_G(t))=1$ for all $t \in D^*$ by Proposition
\ref{restriction} and  $G=\langle C_G(k)\zs|
\; k \in D_1^*\rangle$ for any four-subgroup $D_1 \leq G$ by \cite{chr2}. \hfill $\Box$

\begin{corollary}
\label{abelian}
$C_G(k)\zs$ is nonabelian for all $k \in D^*$.
\end{corollary}

\proof Assume that $C_G(j)\zs$ is abelian for some $j \in D^*$. Then
$C_G(j)\zs \leq C_G(i)\zs$ since $i \in C_G(j)\zs$. As furthermore $j$ and
$ij$ are conjugate, $G=C_G(i)\zs$ by Proposition \ref{corset}. Contradiction.
\hfill $\Box$   

\begin{corollary}
\label{structure}
Let $j \in D^*$. If\/ $i$ and $j$ are not conjugate in $G$, then
 $C_G(j)\zs \cong \PSL_2(K) \times K^*$.
\end{corollary}

\proof Let $j\in D^*$ be not conjugate to $i$. Then $j$ is conjugate
to $ij$ in $C_G(i)\zs$. Now $C_G(j)\zs$ is not abelian by Corollary
\ref{abelian}. Furthermore $C_G(D)\zs \cong K^* \times K^*$ and there
are only
three possibilities by Proposition \ref{restriction}, namely 
$C_G(j)\zs \cong \GL_2(K)$, $C_G(j)\zs \cong
 \PSL_2(K) \times K^*$ or $C_G(j)\zs \cong \SL_2(K)*\SL_2(K)$. The
first and third case cannot occur since  $i$ and $ij$ are not
conjugate, which proves the claim. \hfill $\Box$ \vspace{1ex}

\noindent We have to distinguish two different cases:
\begin{itemize}
\item[(i)]$C_G(i)$ is connected.
\item[(ii)] $C_G(i)$ is not connected.
\end{itemize}
We are going to show that $G \cong \G_2(K)$ in the first case and $G
\cong \PSp_4(K)$ in the second case. Assume from now on that $D=
\langle i, j \rangle$.

\begin{lemma}
\label{concon}
If\/ $C_G(i)$ is connected, then $G$ contains one conjugacy class of
involutions.
\end{lemma}

\proof $S \leq C_G(i)$, since $i \in Z(S)$. As Sylow 2-subgroups of
$G$ are conjugate, and as $C_G(i)$ contains only two conjugacy classes of
involutions $i$ and $j^{C_G(i)}$ -- all elements of order 4 are conjugate in $\SL_2(K)$ --, it is enough to show that $i$ and $j$ are conjugate. 

Assume that $j$ is not conjugate to $i$ and let $w \in N_G(T)$ such
that $D\langle w \rangle$ is an elementary abelian 2-subgroup of order $8$.
Then $C_G(D)=T\langle w \rangle$ and $w$ inverts $T$. As $i$ and $ij$
are not conjugate in $G$,  $C_G(j)\zs \cong  \PSL_2(K) \times K^*$ by Corollary \ref{structure}.  Let $u \in (N_G(T) \cap
C_G(j)\zs)\backslash T$ be an involution. Then $u \in 
C_G(D)$, i.e.\  since $C_G(D) = T\langle w \rangle$, $u=tw$ for some $t \in
T$. Thus $u$ has to invert $T$. Contradiction as $Z(C_G(j)\zs) \leq T$ is
infinite. \hfill $\Box$        

\begin{corollary}
\label{findg}
If\/ $C_G(i)$ is connected, then $G \cong \G_2(K)$.
\end{corollary}
 
\proof If $C_G(i)$ is connected, then $G$ contains one conjugacy class of
involutions by Lemma \ref{concon}. Since furthermore the centralisers of
involutions in $\G_2(K)$ are isomorphic to $\SL_2(K) * \SL_2(K)$, the claim
follows by Theorem \ref{maing2}. \hfill $\Box$

\begin{lemma}
\label{normandall}
 $N_G(L_1) = C_G(i)\zs C_G(C_G(i)\zs)$.
\end{lemma}

\proof As $C_G(i)\zs C_G(C_G(i)\zs) \leq N_G(L_1)$, we only need to
prove the reverse inclusion. Let $g \in N_G(L_1)$. As $Z(L_1)=\langle
i \rangle$, $g \in C_G(i)$ and $g \in N_G(L_2)$ as well by \cite[7.1]{BN}.   
By Lemma \ref{ratherin} and Lemma \ref{cicacen}, $g \in C_G(L_1)L_1 \cap
C_G(L_2)L_2$.  Thus there exists $l_n \in L_n$ for $n=1,2$ such that $ gl_n
\in C_G(L_n)$.  As $L_1$ and $L_2$ commute, $gl_1l_2 \in C_G(L_1) \cap
C_G(L_2)$ and hence $g \in C_G(C_G(i)\zs)C_G(i)\zs$. \hfill $\Box$ 

\noindent\\Set $K_s:=C_G(C_G(s)\zs)$ for all involutions $s \in G$.

\begin{lemma}
\label{control}
Let $s \in I(G)$. Then  $K_s \cap C_G(s)\zs = Z(C_G(s)\zs)$. Furthermore
\begin{itemize}
\item[(i)] $K_i$ is a finite group,
\item[(ii)]$I(K_i)=i$,
\item[(iii)] $K_i \cap K_s =1$ for all $s \in D^*$ such that $i \neq s$.
\end{itemize}
\end{lemma}

\proof As $Z(C_G(i)\zs)=\langle i \rangle$, $K_i \leq C_G(i)$ and $K_i
\cap C_G(i)\zs = \langle i \rangle$. Hence $K_i\zs \leq
C_G(i)\zs \cap K_i = \langle i \rangle$ and $K_i$ is a finite
group, proving $(i)$.
 
To show $(ii)$ let $t \in I(K)$ be an involution. Then $O(C_G(t))=1$ by
\cite{chr2}. Assume that $C_G(t)\zs > C_G(i)\zs$. Then
$C_G(t)\zs \cong \G_2$ or $C_G(t)\zs \cong \PSp_4(K)$ by \cite{nato} as $pr(C_G(t)\zs)=pr(G)=2$. However this would imply
that $C_G(j)\zs, C_G(ij)\zs \leq C_G(t)\zs$ by Proposition \ref{restriction} and
$G=C_G(t)\zs$ by
Proposition \ref{corset}. Contradiction. Thus $C_G(t)\zs = C_G(i)\zs =
C_G(it)\zs$ and $t=i$, as otherwise $\langle i, t \rangle$ is a
four-subgroup and $G=C_G(i)\zs$ by Proposition \ref{corset}.  

Let finally $s \in D^*$ such that $s \neq
i$. As $C_G(i)\zs, C_G(s)\zs \leq C_G(K_i \cap K_s)$ and $si$ is
conjugate to $s$ in $C_G(i)\zs$, \cite[5.14]{nato} and Proposition
\ref{corset} imply that $C_G(K_s \cap K_t)\zs =G$. As $G$ is simple, 
$K_i \cap K_s=1$.  \hfill $\Box$

\begin{corollary}
\label{precision}
$C_G(i) = (C_G(i)\zs * K_i) \semi \langle v \rangle$ for any $v \in G$ such
that $L_1^v=L_2$. In this case $v^2 \in C_G(i)\zs K_i$ and we can choose $v
\in N_G(T)$ to be a 2-element.
\end{corollary} 

\proof Assume that there
 exists an element $v \in G$ such that $L_1^v=L_2$. As
$Z(L_1)=Z(L_2)=\langle i \rangle$,  $v \in
 C_G(i)\backslash C_G(i)\zs$, $L_2^v=L_1$ by \cite[7.1]{BN} and $v^2
  \in C_G(i)\zs K_i$ by Lemma \ref{normandall}. 

Let $x \in C_G(i)$. There are two possibilities: Either $x \in
 N_G(L_1) = C_G(i)\zs * K_i$ by Lemma \ref{normandall} and Lemma
 \ref{control} or $L_1^x=L_2$ and $L_2^x=L_1$ by \cite[7.1]{BN}.
 In the second case $xv  \in N_G(L_1) =
 C_G(i)\zs K_i$ by Lemma \ref{normandall} again. Thus $x \in
 (C_G(i)\zs * K_i) \semi \langle v \rangle$ in all cases.

Let finally $d(v)=V \times F$, where $V$ is a connected divisible group and $F
= \langle f \rangle$ a finite cylic group by \cite[ex.\ 10, p.\
93]{BN}. Then $V \leq C_G(i)\zs$ and $f^2 \in N_G(L_1)$. Thus $f$ has
even order. Let $o(f)=2^km$ where $k \geq 1$ and $m \in \nn$ is
odd. Then $f^m$ is a 2-element such that $L_1^{f^m}=L_2$. Furthermore
$f^m$ is contained in a Sylow 2-subgroup of $C_G(i)$ and as they are
all conjugate to each other, we may assume that $f^m \in N_G(T)$.
\hfill $\Box$  

\begin{lemma}
\label{wemove}
If all involutions in $D$ are conjugate, then $K_i =
\langle i \rangle$. 
\end{lemma}

\proof Assume that all involutions in $D$ are conjugate. Then there exists an
element $y \in N_G(T)$, such that $i^y =j$ by \cite[10.22]{BN}.  Let
$w \in N_G(T) \cap C_G(i)\zs$ such that $E:=D\langle w \rangle$ is an
elementary abelian 2-subgroup of order $8$. $w$ inverts $T$ and all
other involutions that invert $T$
are contained in $wC_G(T)=wTK_i$ by Corollary
\ref{precision}. However, $I(wTK_i)=wT$ by Lemma \ref{control} and
$w^y \in wT$. $K_i$ acts hence on $C_G(j)\zs$ centralising
$T\langle w \rangle$. Thus $K_i \leq C_G(j)\zs K_j$ by Corollary
\ref{precision}. Let $k \in K_i$ and $x \in
C_G(j)\zs$, $k_1 \in K_j$ such that $k=xk_1$. Then $x=kk_1^{-1} \in
C_G(E) \cap C_G(j)\zs =E$. Thus $K_j \leq EK_i$. Furthermore
$K_{ij}=K_j^{w_1} \leq EK_i$, where $w_1 \in L_1 \cap N_G(E)$ is an
involution conjugating $j$ and $ij$. Let now $k:=|K_i|$. As all involutions in
$D$ are conjugate, $k=|K_l|$ for all $l \in D^*$. Furthermore $|EK_i|=4k$ as
$E \cap K_i = \langle i \rangle$. On the other hand $K_i, K_j, K_{ij} \leq
EK_i$ and all three subgroups normalise each other. Hence $k^3 \leq 4k$ by
Lemma \ref{control} and $k=2$. Thus $K=\langle i \rangle$, if all involutions
in $D$ are conjugate. \hfill $\Box$

\begin{proposition} 
\label{exclusion}
If all involutions in $D$ are conjugate, then
$C_G(i)$ is connected. 
\end{proposition}

\proof Suppose that all involutions in $D$ are conjugate. Then
$K_i=\langle i \rangle$ by Lemma \ref{wemove}. Assume that $C_G(i)$ is not
connected. Then there exists a 2-element $v \in
N_G(T)$ such that $C_G(i) = C_G(i)\zs \semi \langle v
\rangle$ where $L_1^v=L_2$ by Lemma \ref{precision}. As $v^2 \in
C_G(i)\zs$ is a 2-element, it is contained in a maximal
torus of $C_G(i)\zs$. As
maximal tori are divisible, there exists an element $t \in C_G(i)\zs$,
such that $t^2=v^2$. Assume that $t=t_1t_2^v$, where $t_1,t_2 \in
L_1$. Then 
\[ (*) \quad 1 = v^2t^{-2} = v^2(t_2^{-2v}t_1^{-2})=vt_2^{-2}vt_1^{-2}.\] 
On the other hand $v^2 \in C_G(v)$ and thus $t_1^2t_2^{2v} =
(t_1^2t_2^{2v})^v = t_2^{2v^2}t_1^{2v}$. This implies that
\[t_2^{-2v^2}t_1^2=t_1^{2v}t_2^{-2v} \in L_1 \cap L_2 = \langle i
\rangle\] 
Hence either $t_1^2=t_2^2$ or $t_1^2=it_2^2$. Set $s:=t_1^2$. If
$t^2=ss^v$, then $(*)$ implies that
$vs^{-1}$ is an involution. If $t^2=iss^v$, then $(vs^{-1})^2=i$. By
replacing $v$ with a conjugate of $vs^{-1}$, we may assume that $v \in
N_G(T)$, such that $v^2 \in \langle i \rangle$.  

As all involutions in $D$ are conjugate
and $j$ is conjugate to $ij$ in $C_G(i)$, there exists an
element $y \in N_G(T)$, such that $i^{y^2} =j^{y} =ij$ by
\cite[10.22]{BN}, where $y^3 \in T$. Furthermore $v \in C_G(D)$ and
$v$ does not invert $T$. This implies that $v
\in C_G(j)\backslash C_G(j)\zs$ and $v \in C_G(ij) \backslash
C_G(ij)\zs$. Thus there exist elements $x_1, x_2 \in C_G(D)$ such that
$v^y=vx_1$ and $v^{y^2}=vx_2$, where $x_1 \in C_G(j)\zs$ and $x_2 \in
C_G(ij)\zs$ by Corollary \ref{precision}. 

Let $P:=C_T(v)\zs$. Then $P=\{ll^v| l \in T \cap L_1\}$ and
$C_T(v)=P\langle i \rangle$. On the other hand $x_n$ either inverts
$T$ or centralizes it for some $n=1,2$. If $x_n$ would centralize $T$ for
$n=1,2$, then $P=C_T(vx_n)\zs=(C_T(v)^{y^n})\zs =
P^{y^n}$. Contradiction as $P$ contains a unique involution from
$D$. Thus $x_1$ and $x_2$ invert $T$.   

Let $T[v]:=\{t \in T| t^v=t^{-1}\}$. Then $T[v] \leq T$ and
$T[v]\zs=\{ll^{-v}| l \in T \cap L_1\}$. Obviously $T=PT[v]$. As $x_n$
inverts $T$ for $n=1,2$, $C_T(x^{y^n})=C_T(vx_n)=T[v]$. Thus $T[v]\zs
= (C_T(v^y)\zs)^y = (T[v]\zs)^y$. Contradiction again, as $T[v]\zs$
contains a unique involution from $D$. Hence
$C_G(i)=C_G(i)\zs$ is connected.  \hfill $\Box$  

\begin{proposition}
\label{connocon}
If\/ $C_G(i)$ is not connected, then $C_G(i)=C_G(i)\zs \semi
\langle u \rangle$, where $u \in C_G(D)$ is an involution such that
$L_1^u=L_2$. 
\end{proposition}

\proof As $C_G(i)$ is not connected, $i$ and $j$ are not conjugate in
$G$ by Proposition \ref{exclusion} and $C:=C_G(j)\zs \cong
\PSL_2(K) \times K^*$ by Corollary \ref{structure}. Let $u \in
(N_G(T) \cap C)\backslash T$ be an
involution. Then $u \in
C_G(D)$ but $u \notin C_G(i)\zs$, as $u$ does not invert
$T$. Thus $C_G(i)=C_G(i)\zs K_i \langle u \rangle$ by Lemma
\ref{precision} and we need to show that $K_i = \langle i \rangle$.

As $K_i$ centralizes $Z(C)$,  $K_i \leq K_jC$ by Lemma \ref{ratherin} and
Lemma \ref{cicacen}. We show that $C_{K_i}(u) = \langle i \rangle$. $u$ acts
on $K_i$. Let $k \in C_{K_i}(u)$. Then $k$ centralizes the elementary abelian 
subgroup $E_1:=\langle  i, j, u \rangle$. As $K_i \leq K_jC$ and
$C_{C}(E_1) =Z(C)\langle u \rangle$, we must have $k \in K_j \langle u
\rangle$ by Lemma \ref{control}.
Thus $k^2 \in K_j \cap K_i =1$ by Lemma \ref{control} and $k \in \langle  i
\rangle $ by Lemma \ref{control} again. Especially $C_{K_i}(u) = \langle i
\rangle$. 

As $|K_i|/|C_{K_i}(u)|=|S|$, where $S =\{[k,u]| \; k \in K_i\}
\subseteq K_i$ consists of elements inverted by $u$ by \cite[ex.\ 17,
p.\ 7]{BN}, $2|S|=|K_i|$. Let $s \in
S$. Then $u^s \in C$ and $u^s=s^uus=us^2$. Furthermore $s$ centralizes
$T$ and thus $u^s \in
N_G(T)\backslash T$ which implies that there exists a $t \in T$ such
that $us^2 = u^s = ut$.  Hence $s^2=t \in K_i \cap T = \langle i
\rangle$. Assume that
$s^2=i$. Then $u^s=ui$. As $\langle j, u \rangle$ is conjugate to $D$
in $C$, $u$ and $uj$ belong to different conjugacy classes. However, there are
at most two conjugacy classes of involutions in $uC_G(i)\zs$ by elementary computation, namely  $u^T$ and $(ui)^T$. Contradiction. Hence $s^2=1$ by Lemma
\ref{control} and $s \in \langle i
\rangle$. Then $|K_i| = 2|S| \leq 4$.  

We finally prove that $|S|=1$. Assume that $|S|=2$. Then $|K_i|=4$ and
there exists an element $k \in K_i$ of order $4$ by Lemma
\ref{control}. As $C_{K_i}(u)=\langle i \rangle$,
$k^u=k^3=k^{-1}$. Especially $u^k = k^uuk = uk^2 = ui$. Contradiction
as before and we are done. \hfill $\Box$    

\begin{proposition}
\label{connocen}
If\/ $C_G(i)$ is not connected, $C_G(j)\cong (\GL_2(K)/\langle
-I \rangle) \semi
\langle w \rangle$ where $-I = \left( \begin{smallmatrix} -1
     & 0 \\ 0 & -1 \end{smallmatrix}\right)$ and $w \in C_G(D) \cap
 C_G(i)\zs$ is an involution that acts as an inverse-transpose
automorphism on $\GL_2(K)$, i.e.
\[ \lemt a & b \\ c & d \rim ^{ w} = (ad -bc)^{-1}\lemt d & - c \\ - b & a
\rim \]
 for any $\left( \begin{smallmatrix} a & b \\ c & d
   \end{smallmatrix}\right) \in \GL_2(K)$.
\end{proposition}

\proof As $C_G(i)$ is not connected, $C_G(j)\zs \cong \PSL_2(K)
\times K^* \cong \GL_2(K)/\langle -I
\rangle$ by Proposition \ref{exclusion} and Corollary
\ref{structure}. Set $C:=C_G(j)\zs$.
 
Let $w \in C_G(D) \cap C_G(i)\zs$ be an involution that inverts
$T$. We show that $C_G(j)=C  \semi \langle w \rangle$. As $w$ inverts
$T$, $w \notin C$
and $C_G(j) \geq C \semi \langle w \rangle$. To prove the other
inclusion let $\ol C:=C/Z(C)$. By Lemma \ref{ratherin}
$C_G(j)=C_{C_G(j)}(\ol C)C$. Let $h \in
C_{C_G(j)}(\ol C)$. Then $T^h \leq TZ(C)=T$ and $h \in N_G(T)$. As $h
\in C_G(j)$ and $i$ and $ij$ are not conjugate, $h \in
C_G(D) = T \langle u, w \rangle$. Thus, as $T \langle u \rangle \leq
C$, $C_G(j)=C\langle w\rangle = C_{C_G(j)}(\ol C)C$.      

As $w$ inverts $T$, $w \notin C_{C_G(j)}(\ol C)$. Thus there exists $v
\in C_G(D)$ such that $wv \in C_{C_G(j)}(\ol C)$. $w$
inverts $T$ and $v \in C$, hence $wv$ inverts $Z(C)$. Furthermore
$[wv,C] \subseteq Z(C)$ and $C=Z(C)*C_C(wv)$ by \cite[ex.\ 10, p.\
98]{BN} where $Z(C) \cap C_C(wv)= \langle j \rangle$. Set
$H:=C_C(wv)\zs$. As $C \cong \PSL_2(K) \times K^*$, $H \cong
PSL_2(K)$ and $C=Z(C) \times H$. We may assume that $v \in H$. 
 
Let $x \in C$. Then there exists $z \in Z(C)$ and $h \in H$ such that
$x=zh$. Furthermore $x^{wv}=z^{-1}h$, and $x^w=z^{-1}h^v$. Thus
$h^w=h^v$ for all $h \in H$. Set $T_1=H \cap T$. As $w$ inverts $T_1$,
$v$ has to invert $T_1$ as well. Hence $v \in I(N_G(T_1)\backslash
T_1)$ and after maybe conjugating $w$ with an element from $T$, we may
assume that $v$ corresponds to $\left( \begin{smallmatrix} 0 & 1 \\
    -1 & 0  \end{smallmatrix}\right)$. This implies the
claim. (Compare Section 3.2.) \hfill $\Box$ \vspace{1ex}

{\it Proof of Theorem \ref{super}.}\/ Let $G$ be as in Theorem
\ref{super}. If $C_G(i)$ is connected, then $G \cong \G_2(K)$ by
Corollary \ref{findg}. If on the other hand $C_G(i)$ is not connected,
$C_G(i)$ and $C_G(j)$ are isomorphic to centralizers of involutions in
$\PSp_4(K)$ by Proposition \ref{connocon} and Proposition
\ref{connocen}. Thus $G \cong \PSp_4(K)$ by Theorem \ref{main}.  \hfill $\Box$

\begin{center} \bf ACKNOWLEDGMENTS \end{center}
This paper is part of the author's Ph.D. thesis. She would like to
thank the Landesgraduiertenf\"orderung Baden-W\"urttemberg for their
financial support. She is very grateful to Alexandre V. Borovik, her
supervisor during her stay in Manchester, who gave all the help one
could wish for and without whom this work would not exist.


\begin{thebibliography}{99}
 \bibitem{alt94} T.~Alt{\i}nel, Groups of Finite Morley Rank with
   Strongly Embedded Subgroups, {\it 
 J.\ Algebra} {\bf 180} (1996), 778-807.   


\bibitem{mixed} T.~Alt{\i}nel / A.~Borovik / G.~Cherlin, Groups of
  mixed type, {\it J.\ Algebra} {\bf 192} (1997), 524-571. 
 
  \bibitem{central} T.~Alt{\i}nel / G.~Cherlin, On central extensions
  of algebraic groups, to appear in J. Symb. Logic.  

 \bibitem{chr1} C.~Altseimer, A Characterisation of $\PSp(4,K)$, to
 appear in Comm. Algebra.

\bibitem{chr2} C.~Altseimer, Strongly Embedded Subgroups of Groups of
  Odd Type, in preparation.

\bibitem{diss} C.~Altseimer, ``Contributions to the Classification of Tame K*-groups of Odd Type and Other Applications of 2-local Theory'', Ph.D. thesis, Freiburg 1998.

\bibitem{qua} C.~Altseimer / A.~Berkman, Quasi- and pseudounipotent groups
of \fmr, Preprint of the Manchester Centre for Pure Mathematics 1998/5.

\bibitem{rusbor}  A.~Borovik, 3-local Characterization of
  $G_2(2^n)$, preprint of the Akademy NAUK CCCP, Institut
  Mathematiki, Nowosibirsk, 1980, (in Russian). 


\bibitem{nato} A.~Borovik,  Simply locally finite groups of 
 finite Morley Rank and odd type, {\it in} ``Proceedings of NATO ASI on 
 Finite and Locally Finite Groups'',  Istanbul, 1994, 248-284.  

 \bibitem{BN} A.~Borovik / A.~Nesin, ``Groups of 
 Finite Morley Rank'',  Oxford University Press, 1994.  

\bibitem{hum75} J.~Humphreys, ``Linear Algebraic Groups'',  
 Springer-Verlag, New York Inc., 1975. 

\bibitem{kramer} L.~Kramer / K.~Tent / H.~Van Maldeghem, Simple Groups of Finite Morley Rank and Tits Buildings, to appear in Israel J. Math.

\end{thebibliography}
\end{document}